\documentclass{article}
\begin{document}
\centerline{\textbf{\normalsize Evaluations of Derivatives of Jacobi Theta Functions in the origin}}\vskip .10in

\[
\]
\centerline{\bf Nikos Bagis}

\centerline{Department of Informatics}

\centerline{Aristotele University of Thessaloniki}
\centerline{Thessaloniki, Greece}
\centerline{nikosbagkis@hotmail.gr}

\[
\]

\begin{quote}

\centerline{\bf abstract\rm}   
In this article using Ramanujan's theory of Eisenstein series we evaluate completely the derivatives of the theta functions $\vartheta_1^{(2\nu+1)}(z)$ and $\vartheta_4^{(2\nu)}(z)$ in the origin in closed polynomials forms using only the first three Eisenstein series of weights 2,4, and 6. 

\[
\]
     
\bf keywords: \rm{theta functions; Ramanujan; derivatives; evaluations;}

\end{quote}

\section{Introduction}
Following Berndt in [3], let $q=e^{2\pi i\tau}$, where $\tau$ is in the upper half-plane $H$, and write 
\begin{equation}
\Phi_{\nu}(q):=\sum^{\infty}_{n=1}\frac{n^{\nu}q^n}{1-q^n}
\end{equation}
and
\begin{equation}
E_2(\tau)=1-24\Phi_{1}(q)-\frac{3}{\pi y}
\end{equation}
\begin{equation}
E_{\nu}(\tau)=1-\frac{2\nu}{B_{\nu}}\Phi_{\nu-1}(q)
\end{equation}
for $n>2$ even and where $y=Im(\tau)>0$, $B_n$ denotes the $n$th Bernoulli number (see [2]).\\ 
Ramanujan then defines
\begin{equation}
P=P(q):=1-24\sum^{\infty}_{n=1}\frac{nq^n}{1-q^n}
\end{equation}
\begin{equation}
Q=Q(q):=1+240\sum^{\infty}_{n=1}\frac{n^{3}q^n}{1-q^n}
\end{equation}
and
\begin{equation}
R=R(q):=1-504\sum^{\infty}_{n=1}\frac{n^{5}q^n}{1-q^n}
\end{equation}
According to Ramanujan we can evaluate (see [3] pg. 319 ) every $\Phi_{2\nu+1}(q)$ in polynomials of $Q$ and $R$.\\
Also we will need the Dedekind's eta-function which is 
\begin{equation}
\eta(\tau):=q^{1/24}\prod^{\infty}_{n=1}(1-q^n)
\end{equation}  
In the literature (see [7],[1]) the first elliptic theta function   
is defined by 
\begin{equation}
\vartheta_1\left(z,t\right):=\sum^{\infty}_{n=-\infty}(-1)^{n-1/2}q_1^{(n+1/2)^2}e^{(2n+1)iz}
\end{equation}
\begin{equation}
\vartheta_4(z,t)=\sum^{\infty}_{n=-\infty}(-1)^nq_1^{n^2}e^{2niz}
\end{equation}
where $t=2\tau$ and hence $q=e^{2\pi i\tau}=q_1=e^{\pi i t}$.
Also hold the following relations
\begin{equation}
\vartheta_1\left(z,t\right)=2\sum^{\infty}_{n=0}(-1)^nq_1^{(n+1/2)^2}\sin[(2n+1)z]
\end{equation}
\begin{equation}
\vartheta_4(z,q)=1+2\sum^{\infty}_{n=1}(-1)^nq_1^{n^2}\cos(2nz)
\end{equation}
Since our attention will be focused on $z$ and we do not make extended use of $q$ or $q_1$ the reader may check everytime he wants, the correspondence, numerically with the program he use.\\
Hence we work with the notations $$\vartheta_1(z)=\vartheta_1\left(z,\frac{t}{2}\right)=\vartheta_1(z,\tau)$$
$$
\vartheta_4(z)=\vartheta_4\left(z,\frac{t}{2}\right)=\vartheta_4(z,\tau)
$$
The sine Fourier series of the first derivative of $\vartheta_1(z)$ is:\\
If $|q|<1$, then  
\begin{equation}
\frac{\vartheta_1'(z)}{\vartheta_1(z)}=\cot(z)+4\sum^{\infty}_{n=1}\frac{q^{2n}}{1-q^{2n}}\sin(2nz) .
\end{equation}     
Here we take $Im(z)=0$.
\[
\]
\textbf{Theorem 1.} (Ramanujan) (see [6])\\
Let
\begin{equation}
S_3(m)=\sum_{a\equiv1(mod 4)}q^{\alpha^2/8}\alpha^m
\end{equation}
Then
\begin{equation}
S_3(1)=\eta^3(\tau)
\end{equation}
\begin{equation}
S_3(3)=\eta^3(\tau)P
\end{equation}
\begin{equation}
S_3(5)=\eta^3(\tau)\frac{(5P^2-2Q)}{3}
\end{equation}
\begin{equation}
S_3(7)=\eta^3(\tau)\frac{35P^3-42PQ+16R}{9}
\end{equation}
and in general
\begin{equation}
S_3(2m+1)=\eta^3(\tau)\sum_{i+2j+3k=m}b_{ijk}P^iQ^jR^k
\end{equation}
\[
\] 
In this article we will find a general way to evaluate $b_{ijk}$ and concecuently the derivatives $\vartheta^{(2m+1)}_1(z)$ in the origin using only $P$, $Q$ and $R$. The same thing we do and with $\vartheta_4^{(2\nu)}(z)$. For to complete our purpose we use the following Theorem of Ramanujan for the evaluation of Eisenstein series:
\[
\]
\textbf{Theorem 2.} (Ramanujan) (see [3] chapter [15])\\    
Let $n$ be integer greater of 1, then  
\begin{equation}
S_{1,{2n}}=\frac{(-1)^{n-1}B_{2n}}{4n}E_{2n}(\tau)
\end{equation}
Also if $n$ is even and exceeding 4, then (see and relation (3))
$$
-\frac{(n+2)(n+3)}{2n(n-1)}S_{1,{n+2}}=-20\left(^{n-2}_{2}\right)S_{1,4}S_{1,{n-2}}+$$ 
\begin{equation}
\sum^{\left[(n-2)/4\right]}_{k=1}{}^{'}\left(^{n-2}_{2k}\right)\left[(n+3-5k)(n-8-5k)-5(k-2)(k+3)\right]S_{1,{2k+2}}S_{1,{n-2k}}
\end{equation}
The prime on the summation sign indicates that if $(n-2)/4$ is an integer, then the last term of the sum is to be multiplied by $\frac{1}{2}$.

\section{The $\vartheta^{(2\nu+1)}_1(0)$ derivatives}

\textbf{Lemma 1.}\\
For $\nu$ non negative integer
\begin{equation}
\vartheta_1^{(2\nu+1)}(0,\tau/2)=2(-1)^{\nu}S_3(2\nu+1)
\end{equation}
\textbf{Proof}\\
Recall the Ramanujans Theorem 1 and differentiate (9) to write
$$ 
\vartheta^{(2\nu+1)}_1(0,\tau/2)=2(-1)^{\nu}\sum^{\infty}_{m=0}(-1)^mq^{(m+1/2)^2/2}(2m+1)^{2\nu+1}=
$$
$$
=2 (-1)^{(\nu)}\sum^{\infty}_{m=-\infty}q^{(4m+1)^2/8}(4m+1)^{2\nu+1}=2(-1)^{\nu}S_3(2\nu+1)=
$$
$$
=2(-1)^{\nu}\eta^3(\tau)\sum_{i+2j+3k=\nu}b_{ijk}P^iQ^jR^k
$$
and the proof is complete.
\[
\]
\textbf{Lemma 2.}
\begin{equation}
\frac{d^{2\nu}}{dz^{2\nu}}\left(\sin(2nz)\sin((2m+1)z)\right)_{z=0}=\frac{(-1)^{\nu}}{2}\left[(2m+1-2n)^{2\nu}-(2m+1-2n)^{2\nu}\right]
\end{equation} 
\textbf{Proof.}\\
It is
\begin{equation}
\sin(z)=\sum^{\infty}_{n=0}\frac{(-1)^{n}z^{2n+1}}{(2n+1)!}
\end{equation}  
Observe first that 
\begin{equation}
\sigma_1(c,n)=\frac{d^{\nu}}{dz^{\nu}}\left(\sin(cz)z^n\right)_{z=0}=n!\cdot i^{-n-1}\frac{1-(-1)^n}{2}\left(^{\nu}_n\right)
\end{equation}
The result follows from (23) and (24). 
\[
\]
\textbf{Lemma 3.}
$$
\frac{d^{\nu}}{dz^{\nu}}\left(\cot(z)\sin(cz)\right)_{z=0}=
$$
\begin{equation}
=\frac{i^{\nu}(1+(-1)^{\nu})}{2}\sum^{\nu-1}_{n=1}\frac{\sigma_1(c,n)}{n!}\left(\frac{d^n}{dz^n}(\cot(z)-1/z)\right)_{z=0}+\frac{i^{\nu}(1+(-1)^{\nu})}{2}\frac{c^{\nu+1}}{\nu+1}
\end{equation} 
\textbf{Proof.}\\
It is
\begin{equation}
\cot(z)=\frac{1}{z}-\sum^{\infty}_{n=1}\frac{(-1)^{n-1}2^{2n}B_{2n}}{(2n)!}z^{2n-1}
\end{equation}  
The result follows as in Lemma 2 using (24) and (26). 
\[
\]
\textbf{Theorem 3.}\\
For $\nu$ not negative integer holds
$$
S_3(2\nu+1)=-\sum^{2\nu-1}_{n=1}\frac{\left|B_{n+1}\right|}{n+1}2^{n+1}\left(^{2\nu}_n\right)S_3(2\nu-n)+\frac{1}{2\nu+1}S_3(2\nu+1)-
$$
\begin{equation}
-2\sum^{\nu}_{l=1}\left(^{2\nu}_{2l-1}\right)2^{2l}S_3(2(\nu-l)+1)\Phi_{2l-1}(q)
\end{equation}
Moreover we have 
\begin{equation}
\vartheta^{(2\nu+1)}_1(0)=2(-1)^{\nu}S_3(2\nu+1)
\end{equation}
\textbf{Proof.}\\
Form (10) we have
$$
\vartheta_1'(z)=\vartheta_1(z)\cot(z)+4\sum^{\infty}_{n=1}\frac{q^{2n}}{1-q^{2n}}\sin(2nz)\vartheta_1(z)
$$
Expanding the $\vartheta_1(z)$ in series (10) and taking the $2\nu$ derivative with respect to $z$ in zero we have
$$
\vartheta^{(2\nu+1)}_1(0)=2\sum^{\infty}_{m=0}(-1)^mq^{(m+1/2)^2}\left(\cot(z)\sin((2m+1)z)\right)^{(2\nu)}_{z=0}+
$$
$$
+8\sum^{\infty}_{n=1}\frac{q^{2n}}{1-q^{2n}}\sum^{\infty}_{m=0}(-1)^mq^{(m+1/2)^2}\left(\sin(2zn)\sin((2m+1)z)\right)^{(2\nu)}_{z=0}
$$
Write
$$
A(\tau)=2\sum^{\infty}_{m=0}(-1)^mq^{(m+1/2)^2}\left(\cot(z)\sin((2m+1)z)\right)^{(2\nu)}_{z=0}
$$
and
$$
B(\tau)=8\sum^{\infty}_{n=1}\frac{q^{2n}}{1-q^{2n}}\sum^{\infty}_{m=0}(-1)^mq^{(m+1/2)^2}\left(\sin(2zn)\sin((2m+1)z)\right)^{(2\nu)}_{z=0}
$$
Then from Lemma 3 we have
$$
A(\tau/2)=
$$
$$=(-1)^{\nu}\sum^{2\nu-1}_{n=1}\frac{d^n}{dx^n}\left(\cot(z)-1/z\right)_{z=0}i^{-n-1}\left(\frac{1-(-1)^n}{2}\right)\left(^{2\nu}_{n}\right)S_3(2\nu-n)+$$
$$
+\frac{(-1)^{\nu}}{2\nu+1}S_2(2\nu+1)
$$
From Lemma 2 we have
$$
B(\tau/2)=(-1)^{\nu-1}\sum^{\nu}_{l=1}\left(^{2\nu}_{2l-1}\right)2^{2l}S_3(2(\nu-l)+1)\sum^{\infty}_{n=1}\frac{q^nn^{2l-1}}{1-q^n}
$$
Using the above identities we get the result.
\[
\] 
Using Theorems 2 and 3 the relations (15),(16) and (17) follow from (14). We give also as examples evaluations for higher order of $\nu$: 

\section{Evaluations}

From Theorem 3 with $\nu=4$ we have to evaluate $\Phi_7(q)=\Phi_7$.
For $\nu=5$ we have the value of $S_3(11)$ but first we have to evaluate $\Phi_9(q)=\Phi_9$. This will be done by using Ramanujan's Theorem 2.\\ 
Here are the first examples of Ramanujan's Theorem 2 
$$
\Phi_7=\frac{1}{480}(Q^2-1)
$$
$$
\Phi_9=\frac{1}{264}(1-QR)
$$
$$
\Phi_{11}=\frac{-691+441Q^3+250R^2}{65520}
$$
$$
\Phi_{13}=\frac{1-Q^2R}{24}
$$
$$
\Phi_{15}=\frac{-3617+1617Q^4+2000QR^2}{16320}
$$
$$
\Phi_{17}=\frac{43867-38367 Q^3 R-5500 R^3}{28728}
$$
...etc
\[
\]
Hence using Theorem 3 we get (where we have set $h=\eta(\tau)$)
\footnotesize
$$
\frac{1}{2}\vartheta_1^{(9)}(0)=S_3(9)=\frac{1}{9} h^3 (35 P^4-84 P^2 Q-12 Q^2+64 P R)
$$
$$
-\frac{1}{2}\vartheta_1^{(11)}(0)=S_3(11)=\frac{1}{9} h^3 (385 P^5-1540 P^3 Q-660 P Q^2+1760 P^2 R+64 Q R)
$$
$$
\frac{1}{2}\vartheta_1^{(13)}(0)=S_3(13)=\frac{1}{27} h^3 (5005 P^6-30030 P^4 Q-25740 P^2 Q^2+552 Q^3+45760 P^3 R+
$$
$$
+4992 P Q R-512 R^2)
$$
$$
-\frac{1}{2}\vartheta_1^{(15)}(0)=S_3(15)=\frac{1}{27} h^3 (25025 P^7-210210 P^5 Q-300300 P^3 Q^2+19320 P Q^3+
$$
$$
+400400 P^4 R+87360 P^2 Q R-3648 Q^2 R-17920 P R^2)
$$
\normalsize
...etc

\section{The $\vartheta^{(2\nu)}_4(0)$ derivatives}

Using the identities (see [7]):
\begin{equation}
\frac{\vartheta^{'}_4(z)}{\vartheta_4(z)}=4 \sum^{\infty}_{n=1}\frac{q^n\sin(2nz)}{1-q^{2n}}
\end{equation}
and the obvious
$$
\frac{q^n}{1-q^{2n}}+\frac{q^{2n}}{1-q^{2n}}=\frac{q^{n}}{1-q^{n}}
$$
we have after derivating $2\nu-1$-times the relation (29):
\begin{equation}
2\cdot4^{\nu}\sum^{\infty}_{n=1}\frac{n^{2\nu-1}q^{2n}}{1-q^{2n}}+\left(\frac{d^{2\nu}\log(\vartheta_4(z))}{dz^{2\nu}}\right)_{z=0}=2\cdot4^{\nu}\sum^{\infty}_{n=1}\frac{q^nn^{2\nu-1}}{1-q^n}
\end{equation} 
But in general holds the following: 
\[
\]
\textbf{Theorem 4.} (Faa di Bruno)\\
Let $f$, $g$ be sufficiently differentiantable functions, then 
\scriptsize
\begin{equation}
\frac{d^{n}f(g(x))}{dx^{n}}
=\sum \frac{n!}{m_1!(1!)^{m_1}m_2!(2!)^{m_2}\ldots m_{n}! (n!)^{m_n}}  f^{(m_1+m_2+\ldots+m_n)}(g(x))\cdot\prod^{n}_{j=1}\left(g^{(j)}(x)\right)^{m_j}
\end{equation}
\normalsize
where the sum is over all n-tuples of nonnegative integers $(m_1,m_2,\ldots,m_n)$ satisfying the constraint 
\begin{equation}
1\cdot m_1+2\cdot m_2+3\cdot m_3+\ldots+n\cdot m_n=n
\end{equation}
\[
\]
Hence from 
\begin{equation}
\frac{d^{2\nu}\log(z)}{dz^{2\nu}}=-\frac{\Gamma(2\nu)}{z^{2\nu}} ,     \nu=1,2,\ldots 
\end{equation}
we get
\begin{equation}
\left(\frac{d^{2\nu}\log(\vartheta_4(z))}{dz^{2\nu}}\right)_{z=0}=-\sum \frac{(2\nu)!}{m_1!(1!)^{m_1}m_2!(2!)^{m_2}\ldots m_{2\nu}!((2\nu)!)^{m_{2\nu}}}\times 
$$
$$
\times\frac{\Gamma\left[2(m_1+m_2+\ldots+m_{2\nu})\right]}{[\vartheta_4(0)]^{2(m_1+m_2+\ldots+m_{2\nu})}}\cdot\prod^{2\nu}_{j=1}\left(\vartheta^{(j)}_4(0)\right)^{m_j}
\end{equation}
\[
\]
Let $x$ be the elliptic modulus and $\bf K\rm$ the complete elliptic integral of the first kind as defined in [4] pages 101-102 , then if we set 
\begin{equation}
w:=\sqrt{\frac{2}{\pi}(1-x)\bf K\rm}
\end{equation}
we have:
\[
\]
\textbf{Theorem 5.}\\ The derivatives $\vartheta^{(2\nu)}_4(0)$, $\nu=1,2,3,\ldots$ can evaluated from
$$
-\frac{1}{2\cdot 4^{\nu}}\sum \frac{(2\nu)!}{m_1!(1!)^{m_1}m_2!(2!)^{m_2}\ldots m_{2\nu}!((2\nu)!)^{m_{2\nu}}}\frac{\Gamma\left[2(m_1+m_2+\ldots+m_{2\nu})\right]}{[\vartheta_4(0)]^{2(m_1+m_2+\ldots+m_{2\nu})}}\times
$$
\begin{equation}
\times\prod^{2\nu}_{j=1}\left(\vartheta^{(j)}_4(0)\right)^{m_j}
=\frac{B_{2\nu}}{4\nu}(E_{2\nu}(2\tau)-E_{2\nu}(\tau)) 
\end{equation}
\begin{equation}
\vartheta^{(2\nu)}_4(0)=w\sum_{0\leq i,j,k,l,m,n\leq \nu}c_{ijklmn}P^{i}Q^{j}R^{k}P_1^{l}Q_1^{m}R_1^{n}
\end{equation}
where $P_1=P(q)$, $P_2=P(q^2)$, $Q_1=Q(q)$, $Q_2=Q(q^2)$, $R_1=R(q)$, $R_2=R(q^2)$ and $c_{ijk}$ rationals.\\
\textbf{Proof.}\\
i) Relation (36) follows from (34),(30) and (1),(2),(3).\\  
ii) For relation (37) observe that (36) is a recurrence relation both in $\vartheta^{(2\nu)}_4(0)$, $E_{2\nu}(\tau)$ and  $E_{2\nu}(2\tau)$. Hence the $\vartheta^{(2\nu)}(0)$ can evaluated using Theorem 2 from $P(q)$, $Q(q)$, $R(q)$ and $P(q^2)$, $Q(q^2)$, $R(q^2)$ (except for $\nu=1$ which we use (4) instead of (2)). 
\[
\]
\textbf{Examples}

\[
\]

\scriptsize

$$
\vartheta_4^{(2)}(0)=\frac{1}{3}(P_1-P)w
$$ 
$$
\vartheta_4^{(4)}(0)=\frac{1}{15} \left(5 (P-P1)^2-2 Q+2 Q1\right) w
$$ 
$$
\vartheta_4^{(6)}(0)=\frac{1}{63}(-35 P^3+105P^2 P_1-105P P_1^2+35P_1^3+42 P Q-42P_1Q-42 P Q_1+42 P_1 Q_1-16 R+16R_1) w
$$
$$
\vartheta_4^{(8)}(0)=\frac{1}{135} [175 P^4-700P^3 P_1+1050 P^2 P_1^2-700 PP_1^3+175 P_1^4-420 P^2 Q+840 P P_1 Q-420 P_1^2 Q-$$
$$-60 Q^2+420 P^2 Q_1-840 P P_1 Q_1+420P_1^2 Q_1-168 Q Q_1+228 Q_1^2+320 P R-320 P_1 R-320 P R_1+320P_1 R_1]w
$$

\normalsize

\newpage
\[
\]
\centerline{\bf References}\vskip .2in

\noindent

[1] J.V. Armitage W.F. Eberlein, 'Elliptic Functions'. Cambridge University Press (2006).

[2]: M.Abramowitz and I.A.Stegun, 'Handbook of Mathematical Functions'. Dover Publications, New York. 1972.

[3]: B.C.Berndt, Ramanujan`s Notebooks Part II. Springer Verlang, New York (1989)

[4]: B.C.Berndt, Ramanujan's Notebooks Part III. Springer Verlang, New York (1991)

[5]: I.S. Gradshteyn and I.M. Ryzhik, 'Table of Integrals, Series and Products'. Academic Press (1980).

[6]: Heng Huat Chan, Shaun Cooper and Pee Choon Toh, 'Ramanujan's Eisenstein Series and Powers of Dedekind's Eta-Function'. J. London Math. Soc. (2) 75 (2007) 225-242.  

[7]: E.T.Whittaker and G.N.Watson, 'A course on Modern Analysis'. Cambridge U.P. (1927)

\end{document}